\documentclass[15pt, reqno]{amsart}
\usepackage{amsmath}
\usepackage{amssymb}

\numberwithin{equation}{section}

\begin{document}
\setlength{\baselineskip}{1.5em}

{\theoremstyle{plain}
	\newtheorem{theorem}{\bf Theorem}[section]
	\newtheorem{proposition}[theorem]{\bf Proposition}
	\newtheorem{claim}[theorem]{\bf Claim}
	\newtheorem{lemma}[theorem]{\bf Lemma}
	\newtheorem{corollary}[theorem]{\bf Corollary}
}
{\theoremstyle{remark}
	\newtheorem{remark}[theorem]{\bf Remark}
	\newtheorem{example}[theorem]{\bf Example}
}
{\theoremstyle{definition}
	\newtheorem{defn}[theorem]{\bf Definition}
}

%%%%%%%%%%%%%%%%%%%%%%%%%%%%%%%%%%%%%%%%%%%%%%%%%%%%%%%%%%%%%%%%

\newcommand{\spec}{\operatorname{Spec}}
\newcommand{\proj}{\operatorname{Proj}}
\newcommand{\pic}{\operatorname{Pic}}
\newcommand{\Hom}{\operatorname{Hom}}
\newcommand{\ext}{\operatorname{Ext}}
\newcommand{\tor}{\operatorname{Tor}}
\newcommand{\reg}{\operatorname{reg}}
\newcommand{\vol}{\operatorname{Vol}}
\newcommand{\m}{\dashrightarrow}
\newcommand{\smap}{\rightarrow\!\!\!\!\!\rightarrow}
\newcommand{\sfrac}[2]{\frac{\displaystyle #1}{\displaystyle #2}}
\newcommand{\under}[1]{\underline{#1}}
\newcommand{\ov}[1]{\overline{#1}}
\newcommand{\sheaf}[1]{{\mathcal #1}}
\newcommand{\supp}{\operatorname{Supp}}

\def\l{\lambda}
\def\To{\longrightarrow}
\def\height{\operatorname{ht}}
\def\mm{{\frak m}}
\def\pp{{\frak p}}
\def\O{{\mathcal O}}
\def\I{{\mathcal I}}
\def\M{{\mathcal M}}
\def\L{{\mathcal L}}
\def\S{{\mathcal S}}
\def\N{{\mathcal N}}
\def\NN{{\mathbb N}}
\def\PP{{\mathbb P}}
\def\ZZ{{\mathbb Z}}
\def\QQ{{\mathbb Q}}
\def\CC{{\mathbb C}}
\def\T{{\mathcal T}}
\def\H{\tilde{H}}

%%%%%%%%%%%%%%%%%%%%%%%%%%%%%%%%%%%%%%%%%%%%%%%%%%%%%%%%%%%%%%%%%%%%%%%%

\title{Asymptotic behavior of the length of local cohomology}
\author{Steven Dale Cutkosky}\thanks{The first author was partially supported by NSF}
\author{Huy T\`ai H\`a}
\author{Hema Srinivasan}\thanks{The third author was partially supported by NSA}
\author{Emanoil Theodorescu}
\address{Department of Mathematics, University of Missouri, Columbia MO 65201, USA}
\email{cutkoskys@missouri.edu, tai@math.missouri.edu,  srinivasanh@missouri.edu, \newline theodore@math.missouri.edu}

\begin{abstract}
Let $k$ be a field of characteristic 0, $R=k[x_1, \ldots, x_d]$
be a polynomial ring,  and $\mm$ its maximal homogeneous ideal. Let $I \subset R$ be a homogeneous ideal in $R$. In this paper, we show that $$
\lim_{n \to \infty} \sfrac{\l(H^0_\mm(R/I^n))}{n^d}
=\lim_{n \to \infty} \sfrac{\l(\ext^d_R(R/I^n,R(-d)))}{n^d}
$$ 
 always exists. This limit has been shown to be $\sfrac{e(I)}{d!}$  for $m$-primary ideals $I$ in a local Cohen Macaulay ring \cite {ki,ko1,th,th2}, where $e(I)$ denotes the multiplicity of $I$.  
But we find that this limit may not be rational in general.  We give an example for which the limit is an irrational number thereby showing that the lengths of these extention modules may not have polynomial growth.  
\end{abstract}
\maketitle
%%%%%%%%%%%%%%%%%%%%%%%%%%%%%%%%%%%%%%%%%%%%%%%%%%%%%%%%%%%%%%%%%%%%%%%%%

\section*{Introduction}

Let $R = k[x_1, \ldots, x_d]$ be a polynomial ring over a field $k$,
with graded maximal ideal $\mm$,
 and $I \subset R$ a proper homogeneous ideal.  We investigate the asymptotic growth of $\l(\ext^d_R(R/I^n, R))$ as a function of $n$.  When $R$ is a local Gorenstein ring and $I$ is an $\mm$-primary ideal, then this is easily seen to be equal to $\l(R/I^n)$ and hence is a polynomial 
in $n$.  A  theorem of Kirby, Kodiyalam and Theodorescu \cite{ki,ko1,th,th2} extends this to $\mm$-primary ideals in 
local Cohen Macaulay rings R. We consider homogeneous ideals in a polynomial ring which 
are not $\mm$-primary and show that a limit exists asymptotically although it can be irrational.  In our setting,  by local duality, 
$$
\l(\ext^d_R(R/I^n, R(-d)))=\l(H^0_\mm(R/I^n))
$$
 and thus this becomes a problem of asymptotic lengths of local cohomology modules. 

In recent years, a great deal 
of interest has been given to investigating asymptotic behavior of algebraic
 invariants of powers of $I$. Cutkosky, Herzog and Trung \cite{cht}, 
and Kodiyalam \cite{ko2} independently proved that $\reg(R/I^n)$ is a
 linear function in $n$ for $n \gg 0$ (see also \cite{ch, ggp}).
 When $I^n$ is replaced by its saturation  $(I^{n})^{\text{sat}}$,
 the problem becomes much subtler. 
It is no longer true that $\reg(R/(I^{n})^{\text{sat}})$ is always asymptotically
 a polynomial in $n$ as shown in  \cite{c}.  
Examples are given in \cite{c,cel}  showing that it is possible for 
$\lim_{n \to \infty} \sfrac{\reg(R/(I^{n})^{\text{sat}})}{n}$ to be an irrational number.
Further,  Cutkosky, Ein and Lazarsfeld \cite{cel} showed that the 
limit $\lim_{n \to \infty} \sfrac{\reg(R/(I^{n})^{\text{sat}})}{n}$
 always exists. Along this theme, Hoa and Hyry \cite{hh} recently
 studied the existence of similar limits where the regularity of $R/I^n$ 
is replaced by its a-invariants.  This paper addresses a closely
 related question. We prove 

\begin{theorem} \label{main}
Let $k$ be a field of characteristic zero, $R = k[x_1, \ldots, x_d]$ be a polynomial ring of dimension $d > 1$, and $\mm$ 
its maximal homogeneous ideal. Let $I \subset R$ be a homogeneous ideal of $R$. Then, the limit
 $$
\lim_{n \to \infty} \sfrac{\l(H^0_\mm(R/I^n))}{n^d}
=\lim_{n \to \infty} \sfrac{\l(\ext^d_R(R/I^n,R(-d)))}{n^d}
$$ 
always exists. 
\end{theorem}

In fact, we prove (in Theorem \ref{main-proof}) that if $R$ is a coordinate ring of a projective
variety which has depth $\ge 2$ at is irrelevant ideal, then the limit
$\lim_{n \to \infty} \sfrac{\l(H^0_\mm(R/I^n))}{n^d}$
exists, and (Corollary \ref{Cor}) if $R$ is Gorenstein, then
$$
\lim_{n \to \infty} \sfrac{\l(H^0_\mm(R/I^n))}{n^d}
=\lim_{n \to \infty} \sfrac{\l(\ext^d_R(R/I^n,R(-d)))}{n^d}.
$$ 

We will also give an example where $\sfrac{\l(H^0_\mm(R/I^n))}{n^d}$ tends to an irrational number
 as $n \to \infty$ (Theorem \ref{irrationality}).
 This, in particular, shows that, just like $\reg(R/(I^{n})^{\text{sat}})$, $\l(H^0_\mm(R/I^n))$ is not asymptotically a polynomial in $n$.

Theorem \ref{main} is proved in Section 1. To do this, we express $\l(H^0_\mm(R/I^n))$ as a sum of two components, the geometric component $\sigma(n)$, and the algebraic component $\tau(n)$, and show that both limits $\lim_{n \to \infty}\sfrac{\sigma(n)}{n^d}$ and $\lim_{n \to \infty}\sfrac{\tau(n)}{n^d}$ exist. For the first limit, we express $\sigma(n)$ as $h^0(Y, \N^n)$ for some line bundle $\N$ over a projective scheme $Y$ of dimension $d$, and investigate the limit $\lim_{n \to \infty} \sfrac{h^0(Y, \N^n)}{n^d}$. For the later one, we write $\tau(n)$ as the Hilbert function of a finitely generated graded $k$-algebra of dimension $(d+1)$. In Theorem \ref{Theorem*}, we use the construction illustrated in \cite{c} by the first author to give an example where the limit proved to exist in Section 1 is an irrational number.

\begin{theorem}\label{Theorem*} There exists a nonsingular projective curve $C\subset
\PP^3_{\CC}$ such that if $I\subset R=\CC[x_1,\ldots,x_4]$ is the defining ideal of $C$,
and $\mm$ is the homogeneous maximal ideal of $R$, then
$$
\lim_{n \to \infty} \sfrac{\l(H^0_\mm(R/I^n))}{n^4} \not\in \QQ. 
$$
\end{theorem}

Of course, Theorem \ref{Theorem*} has a local analog.

\begin{theorem}\label{examp} There exists a regular local ring $S$ of dimension 4 which is essentially
of finite type over the complex numbers ${\CC}$, and an ideal $J\subset S$ such that
$$
\lim_{n \to \infty} \sfrac{\l(\ext^d_S(S/J^n,S)))}{n^d}
$$
is an irrational number. In particular, $\l(\ext^d_S(S/J^n,S)))$ is not a polynomial or
a quasi-polynomial for large $n$.
\end{theorem}

The proofs of Theorems \ref{Theorem*} and \ref{examp} will be given in section 2. In contrast to the example of Theorem \ref{examp}, if $(S,\mm,k)$ is a Cohen-Macaulay  local ring of dimension $d$ and $J$ is an
$\mm$-primary ideal, then  $\l(\ext^d(S/J^n,S)$ is a polynomial of degree $d$ for large $n$ 
(\cite{ki}, \cite{th}). In fact,  
$$
\lim_{n \to \infty} \sfrac{\l(\ext^d_S(S/J^n,S)))}{n^d}=\sfrac{e(I)}{d!}
$$
where $e(I)$ is the multiplicity of $I$ (\cite{th2}).

In the case when $S$ is Gorenstein, this follows easily from
local duality, since
$$
\l(\ext^d(S/J^n,S))=\l(H^0_{\mm S}(S/J^n))=\l(S/J^n)
$$
for all $n$.

Suppose that $I$ is a homogeneous ideal in the coordinate ring $R$ of a projective variety
of depth $\ge 2$ at the irrelevant ideal. Then (c.f. Remark \ref{Remark}), if 
$\height I=d=\text{ dim }(R)$, we have that 
$$
\lim_{n \to \infty} \sfrac{\l(H^0_\mm(R/I^n))}{n^d}=\frac{e(I)}{d!}
$$
where $e(I)$ is the multiplicity of $I$.

In contrast, if $\height I<d$ then
\begin{equation}\label{new}
\lim_{n \to \infty} \sfrac{\l(H^0_\mm(R/I^n))}{n^d}
\end{equation}
does not have such a simple arithmetic interpretation. The example of
Theorem \ref{Theorem*} of this paper is of a height 2 prime ideal in a polynomial ring
$R$ of dimension 4 such that (\ref{new}) is an irrational number.
However, in many cases, such as when $I$ is a regular prime with $\height I<d$
in a polynomial ring $R$ of dimension $d$, we have that the limit (\ref{new})
is 0.

More generally, if $\height I<d$ and the analytic spread $\ell(I)<d$ then
the integral closure $\overline{I^n}$ has no $\mm$-primary component for large $n$
\cite[Theorem 3]{M}.  Since $\text{depth}(R_{\mm})\ge 2$,
$\overline{I^n}=H^0(\text{spec}(R)-\{\mm\},\overline{I^n})$, and
$$
H^0_{\mm}(R/\overline{I^n})\cong H^1_{\mm}(\overline{I^n})\cong H^0(\text{spec}(R)-\{\mm\},
\overline{I^n})/(\overline{I^n})=0
$$
for large $n$. Thus if $I$ is a normal ideal ($I^n=\overline{I^n}$ for all $n$)
 with analytic spread $\ell(I)<d$ we have that $H^0_{\mm}(R/I^n)=0$
for large $n$ and the limit (\ref{new}) is thus 0. In fact, Catalin Ciuperca has shown us that
even if $I$ is not normal, with $\ell(I)<d$, then the limit (\ref{new}) is zero.
%%%%%%%%%%%%%%%%%%%%%%%%%%%%%%%%%%%%%%%%%%%%%%%%%%%%%%%%%%%%%%%%%%%%%%%%%

\section{The existence theorem}

In this section, we prove the main theorem of the paper. We shall start by recalling some notations and terminology, and prove a few preliminary results. 

Suppose $R$ is a graded ring, and $I \subset R$ a homogeneous ideal. The {\it Rees algebra} of $I$ is the subalgebra $R[It]$ of $R[t]$. The Rees algebra $R[It]$ has a natural bi-gradation given by
$$ R[It]_{(m,n)} = (I^n)_m t^n. $$
Suppose $A = \oplus_{m,n \in \ZZ}A_{(m,n)}$ is a bi-graded algebra. For a tuple of positive integers $\Delta = (a,b)$, $A_\Delta = \oplus_{n \in \ZZ}A_{(an,bn)}$ is call a {\it $\Delta$-diagonal subalgebra of $A$}. 

\begin{lemma} \label{finite-gen}
Suppose a domain $R$ is a finitely generated graded $k$-algebra of dimension $\delta$, and $I \subset R$ is a homogeneous ideal generated in degrees $\le d$ such that $\height I \ge 1$. Let $A = R[It]$ be the Rees algebra of $I$ over $R$. \par
{\rm (i)} For any tuple $\Delta = (a,b)$ of positive integers such that $a \ge db$, $A_\Delta$ is a finitely generated graded $k$-algebra. \par
{\rm (ii)} For any tuple $\Delta = (a,b)$ of positive integers such that $a > db$, $\dim A_\Delta = \delta$.
\end{lemma}

\begin{proof} It is easy to see that $A_\Delta = k[(I^b)_a]$ is the $k$-algebra generated by elements of $(I^b)_a$. Thus, (i) is clear. (ii) follows from \cite[Lemma 2.2]{hot} since $R$ is a domain.
\end{proof}

The following Lemma is stated in an example in \cite{la}.

\begin{lemma} \label{limit}
Suppose that $Y$ is a projective variety of dimension $d$ over a field $k$ of characteristic zero, and $\L$ is a line bundle on $Y$. Then,  the limit
$$ \lim_{n \to \infty} \sfrac{h^0(Y, \L^n)}{n^d} $$
exists, and is a positive real number if $\L$ is big.
\end{lemma}

\begin{proof} It follows from  \cite[Theorem 10.2]{I} that $\limsup_{n \to \infty} \sfrac{h^0(Y, \L^n)}{n^d} = 0$ if $\L$ is not big ($\kappa(\L)<d$).
This implies that $\lim_{n \to \infty} \sfrac{h^0(Y, \L^n)}{n^d} = 0$.

 Suppose that $\L$ is big (that is, $\kappa(\L)=d$).
It follows from  \cite[Theorem 10.2]{I} that 
$$
\liminf_{n \to \infty} \sfrac{h^0(Y, \L^n)}{n^d} > 0.
$$
To prove the lemma, it suffices to show that
$$ \limsup_{n \to \infty}\sfrac{h^0(Y, \L^n)}{n^d/d!} = \liminf_{n \to \infty} \sfrac{h^0(Y, \L^n)}{n^d/d!}. $$ 

Let $\varepsilon > 0$ be an arbitrary positive number. By applying the theorem of Fujita \cite{fu} (and from the definition of limsup), there exists a birational morphism $\theta: Z \To Y$ together with an effective $\QQ$-divisor $E$ on $Z$ such that $H = \theta^* \L - E$ is a semiample $\QQ$-divisor with 
$$H^d > \limsup_{n \to \infty}\sfrac{h^0(Y, \L^n)}{n^d/d!} - \varepsilon. $$ 
Let $q$ be the smallest positive integer such that $qE$ is integral (or equivalently, $qH$ is integral). Since $qE$ is effective, there is a natural injection $\O_Z \hookrightarrow \O_Z(qE)$. This gives an injection $\O_Z(lqH) \hookrightarrow \O_Z(lqE + lqH) = \theta^* \L^{lq}$ for any integer $l > 0$. Thus, $h^0(Z, \O_Z(lqH)) \le h^0(Z, \theta^* \L^{lq})$ for any integer $l > 0$. Furthermore, since $qH$ is semiample, there exists $c>0$ such that 
\begin{equation}\label{eq*}
 h^i(Z, \O_Z(lqH)) \le c(l^{d-2})
\end{equation}
 for any $i \ge 1$ (c.f. Corollary 6.7 \cite{F2}). 
Thus 
$$h^0(Z,\O_Z(lqH))= \chi(lqH) +O(l^{d-2})= \sfrac{(lq)^d}{d!}H^d + O(l^{d-1}), $$
where $\chi$ denotes the Euler characteristic.

We have an exact sequence of coherent $\O_Y$ modules
$$
0\rightarrow \O_Y\rightarrow \theta_*\O_Z\rightarrow {\mathcal F}\rightarrow 0
$$
where $\mathcal F$ is supported on a closed subset of $Y$ of dimension $<d$. From the
exact sequences
$$
0\rightarrow {\mathcal L}^{lq}\rightarrow \theta_*\theta^*{\mathcal L}^{lq}\rightarrow {\mathcal F}\otimes{\mathcal L}^{lq}\rightarrow 0
$$
we have
$$
h^0(Y,{\mathcal L}^{lq})=h^0(Z,\theta^*{\mathcal L}^{lq})+O(l^{d-1}).
$$
Hence, 
$$ H^d = \lim_{l \to \infty} \sfrac{h^0(Z, \O_Z(lqH))}{(lq)^d/d!} \le \liminf_{l \to \infty} \sfrac{h^0(Y, \L^{lq})}{(lq)^d/d!}. $$ 
In summary, we have
\begin{align}
\limsup_{n \to \infty}\sfrac{h^0(Y, \L^n)}{n^d/d!} < H^d + \varepsilon \le \liminf_{l \to \infty} \sfrac{h^0(Y, \L^{lq})}{(lq)^d/d!} + \varepsilon. \label{supinf}
\end{align}
Since $\varepsilon$ is taken to be arbitrary, by (\ref{supinf}), to prove the lemma, we only need to show that 
\begin{align}
\liminf_{n \to \infty} \sfrac{h^0(Y, \L^n)}{n^d/d!} = \liminf_{l \to \infty}\sfrac{h^0(Y, \L^{lq})}{(lq)^d/d!}. \label{liminf}
\end{align}

Since $\L$ is big, there exists a constant $n_0$ such that $h^0(Y, \L^n) > 0$ for any $n \ge n_0$ (as follows from \cite[Theorem 10.2]{I}). This implies that, for each $0 \le r < q$, there exists an effective divisor $F_r$ on $Y$ such that $\O_Y(F_r) \cong \L^{n_0+r}$. For $l > 0$, we have exact sequences
$$ 0 \To \L^{lq} \To \L^{lq+n_0+r} \To \O_{F_r} \otimes \L^{lq+n_0+r} \To 0. $$
Taking the long exact sequences of cohomologies, we get
$$ 0 \to H^0(Y, \L^{lq}) \to H^0(Y, \L^{lq+n_0+r}) \to H^0(F_r, \O_{F_r}\otimes \L^{lq+n_0+r}). $$ Thus,
\begin{align} 
\sfrac{h^0(Y, \L^{lq+n_0+r})}{(lq)^d/d!} - \sfrac{h^0(Y, \L^{lq+n_0+r})}{(lq)^d/d!}
\le
\sfrac{h^0(F_r, \O_{F_r}\otimes \L^{lq+n_0+r})}{(lq)^d/d!} . \label{reduction}
\end{align}
Since every component of $F_r$ has dimension $d-1<d$, we have
$$ \lim_{l \to \infty} \sfrac{h^0(F_r, \O_{F_r}\otimes \L^{lq+n_0+r})}{(lq)^d/d!} = 0. $$
(\ref{reduction}) now gives us
$$ 
\liminf_{l \to \infty}\sfrac{h^0(Y, \L^{lq})}{(lq)^d/d!} =
 \liminf_{l \to \infty}\sfrac{h^0(Y, \L^{lq+n_0+r})}{(lq)^d/d!}. 
$$
Moreover, since $n_0$ is fixed and $0 \le r < q$, we have
$$ \liminf_{l \to \infty}\sfrac{h^0(Y, \L^{lq+n_0+r})}{(lq)^d/d!} = \liminf_{l \to \infty}\sfrac{h^0(Y, \L^{lq+n_0+r})}{(lq+n_0+r)^d/d!}. $$
Therefore, 
$$ \liminf_{l \to \infty}\sfrac{h^0(Y, \L^{lq})}{(lq)^d/d!} = \liminf_{l \to \infty}\sfrac{h^0(Y, \L^{lq+n_0+r})}{(lq+n_0+r)^d/d!} $$
 for any $0 \le r < q$. Hence, 
$$ \liminf_{l \to \infty}\sfrac{h^0(Y, \L^{lq})}{(lq)^d/d!} = \liminf_{n \to \infty}\sfrac{h^0(Y, \L^n)}{n^d/d!}. $$
(\ref{liminf}) is proved, and so is the lemma.
\end{proof}

\begin{theorem} \label{main-proof}
Let $R = k[x_1, \ldots, x_p]/J$ be the quotient of a polynomial ring 
$k[x_1,\ldots,x_p]$ over a field $k$ of characteristic zero by a 
homogeneous prime ideal $J$. Let $\mm$ be the maximal homogeneous ideal of $R$.
Suppose that $\text{depth}(R_{\mm})\ge 2$. Let $d$ be the dimension of $R$.

 Let $I \subset R$ be a homogeneous ideal of $R$. Then, the limit $\lim_{n \to \infty} \sfrac{\l(H^0_\mm(R/I^n))}{n^d}$ always exists. 
\end{theorem}

By duality, we obtain the following corollary.

\begin{corollary}\label{Cor} With the notations of Theorem \ref{main-proof}, suppose that $R$
is Gorenstein. Then
$$
\lim_{n \to \infty} \sfrac{\l(H^0_\mm(R/I^n))}{n^d}
=\lim_{n \to \infty} \sfrac{\l(\ext^d_R(R/I^n,R(-d))))}{n^d}
$$ 
\end{corollary}

\noindent {\it Proof of Theorem \ref{main-proof}:}  The problem is trivial if $\height I = 0$.
If $\height I=d$ then there is a very simple proof (see the Remark after this proof). 
Suppose $\height I \ge 1$. Let $\I$ be the ideal sheaf associated to $I$ on $V = \proj R$. From the exact sequence
$$ 0 \To I^n \To R \To R/I^n \To 0, $$
we get $H^0_\mm(R/I^n) = H^1_\mm(I^n)$. The Serre-Grothendieck correspondence gives us the exact sequence
$$ 0 \to I^n \to \oplus_{m \ge 0} H^0(V, \I^n(m)) \to H^1_\mm(I^n) \to 0. $$
$\text{depth}(R_{\mm})\ge 2$ implies that
$$
\bigoplus_{m\ge 0}H^0(V,\I^n(m))=H^0(\text{spec}(R)-\{\mm\},I^n)=(I^n)^*
$$
where $(I^n)^*$ is the intersection of the primary components of $I^n$ which are not
$\mm$-primary.  By the theorem of Swanson \cite{S}
there exists a number $e > 0$ such that $(I^n)_m = (I^n)_m^*$ for any $m \ge en$ and $n \ge 1$. Therefore, we have
\begin{equation}\label{eq3}
 \l(H^1_\mm(I^n)) = \sigma(n) - \tau(n),
\end{equation}
where
\begin{equation}\label{eq2}
\sigma(n) = \sum_{m = 0}^{en} h^0(V, \I^n(m))\text{ and }
\tau(n) = \sum_{m = 0}^{en} \l((I^n)_m).
\end{equation}
We will take $e$ to be bigger than the degrees of homogeneous generators of $I$. The theorem will be proved if we can show that both limits $\lim_{n \to \infty} \sfrac{\sigma(n)}{n^d}$ and $\lim_{n \to \infty}\sfrac{\tau(n)}{n^d}$ exist. 

Let us first consider $\lim_{n \to \infty} \sfrac{\sigma(n)}{n^d}$. Let $\pi: X \To V$ be the blowing up of $V$ along $\I$. Let us denote $\M = \pi^* \O_{V}(1)$ and $\L = \I \O_X$. Let $\eta: Y = \PP(\O_X \oplus \M) \to X$ be the projectivization of the vector bundle $\O_X \oplus \M$ on $X$. Then, $\dim Y = \dim X + 1 = d$. Let $\N = \O_Y(e) \otimes \eta^* \L$. We have
$$ h^0(Y, \N^n) = h^0(X, S^{en}(\O_X \oplus \M) \otimes \L^n) = \sum_{m = 0}^{en} h^0(X, \M^m \otimes \L^n).$$
Furthermore, it follows from \cite[Exercise II.5.9]{H}
 (see also \cite[Lemma 3.3]{cel} and \cite{ht}) that $\pi_* \L^n = \I^n$ for $n \gg 0$. Thus, for $n \gg 0$,
$$ h^0(Y, \N^n) = \sum_{m=0}^{en} h^0(V, \O_{V}(m) \otimes \I^n) = \sum_{m=0}^{en}h^0(V, \I^n(m)) = \sigma(n).$$
By Lemma \ref{limit}, $\lim_{n \to \infty} \sfrac{h^0(Y, \N^n)}{n^d}$ exists. Hence, there exists a limit
$ \lim_{n \to \infty} \sfrac{\sigma(n)}{n^d}.$

Now, let us consider $\lim_{n \to \infty}\sfrac{\tau(n)}{n^d}$. Suppose $I$ is generated by $F_1, \ldots, F_l \in R$ and $\deg F_j = d_j$ for $1 \le j \le l$. Let $S = R[s, F_1t, \ldots, F_lt] \subset R[s,t]$ be the Rees algebra of the ideal $IR[s]$ over the polynomial ring $R[s]$. $S$ can be viewed as a bi-graded ring with $\deg x_i = (1,0)$ for $1 \le i \le d$, $\deg s = (1,0)$, and $\deg F_jt = (d_j,1)$ for $1 \le j \le l$. Take an abitrary element $f \in S$. We observe that $\deg f = (en,n)$ if and only if $f$ has the following form
$$f = \sum_{m_1+\cdots+m_l = n} \sum_{j=0}^{en-d_1m_1 - \cdots - d_lm_l} b_j s^{en-d_1m_1-\cdots-d_lm_l-j}F_1^{m_1} \cdots F_l^{m_l} t^n, $$
where $b_j\in R$ is homogeneous of degree $j$.
Thus the map
$$
\Phi_n:S_{(en,n)}\rightarrow \sum_{m=0}^{en}(I^n)_m
$$
defined by $\Phi_n(f(s,t))=f(0,1)$ is a $k$-vector space isomorphism.
 Hence, 
$$\l(S_{(en,n)}) = \sum_{m=0}^{en}\l((I^n)_m) = \tau(n).$$

Let $T = \oplus_{n=0}^{\infty} S_{(en,n)}$, then $T = S_\Delta$ with $\Delta = (e,1)$. Since $S$ is a finitely generated bi-graded $k$-algebra and $e$ is taken to be bigger than $d_j$ for all $1 \le j \le l$, it follows from Lemma \ref{finite-gen} that $T$ is a finitely generated $k$-algebra and $\dim T = d+1$. Thus, the Hilbert function $H(T,n) = \dim_k T_n = \l(S_{(en,n)})$ is given by a polynomial of degree $d$ in $n$ for $n \gg 0$ with a rational leading coefficient. This implies the existence of the limit
\begin{equation}\label{eq1} \lim_{n \to \infty}\sfrac{\tau(n)}{n^d} \in \mathbb R.\end{equation}
The theorem is proved.
\vskip .2truein

\begin{remark}\label{Remark} Let assumptions be as in the statement of Theorem \ref{main-proof}.
If $\height I=d$ we have that $H^0_{\mm}(R/I^n)=R/I^n$ for all $n$, so that
$\l(H^0_{\mm}(R/I^n))$ is the Hilbert polynomial of $I$ for large $n$. Thus
$$
\lim_{n \to \infty} \sfrac{\l(H^0_\mm(R/I^n))}{n^d}=\frac{e(I)}{d!}
$$
where $e(I)$ is the multiplicity of $I$. In the proof we in fact have that 
$\sigma(n)=\sum_{m=0}^{en}h^0(V,{\O}_V(m))$ in this special case.
 
\end{remark}

%%%%%%%%%%%%%%%%%%%%%%%%%%%%%%%%%%%%%%%%%%%%%%%%%%%%%%%%%%%%%%%%%%%%%%%%%

\section{Irrational asymptotic behaviour}

In this section, we will give  examples, stated in Theorems \ref{Theorem*} and \ref{examp} of the introduction,
in which the limit proved to exist in Section 1 is an irrational number. This exhibits how complicated the length $\l(H^0_{\mm}(R/I^n))$ can be asymptotically. In fact, we will show that the construction given by the first author in \cite{c} provides an example. 

Let $S$ be a K3 surface defined over the complex field $\CC$ with $\pic(S) \cong \ZZ^3$. We can therefore identify $\pic(S)$ with integral points $(x,y,z) \in \ZZ^3$. Take $S$ to be the K3 surface which has the intersection form 
\begin{equation}\label{eq4}
q(x,y,z) = 4x^2 - 4y^2 - 4z^2,
\end{equation}
 where $q(D) = D^2$ for any divisor $D \in \pic(S)$. Such a surface $S$ exists
 as shown in \cite{c}. It is shown there that a divisor $D$ on $S$
 is ample if and only if it is in the
interior of
$$
\overline{NE}(S)=\{(x,y,z)\in {\mathbb R}^3\mid 4x^2-4y^2-4z^2\ge 0, x\ge 0\}.
$$
Moreover, $S$ is embedded into $\PP^3$ by the divisor $H=(1,0,0)$. Suppose $(a,b,c) \in \ZZ^3$ is such that 
\begin{align} 
\left\{ \begin{array}{rcl} a & > & 0, \\
a^2 - b^2 - c^2 & > & 0, \\
\sqrt{b^2+c^2} & \not\in & \QQ. \end{array} \right. \label{abc}
\end{align} 
Since $(a,b,c)$ is in the interior of $\overline{NE}(S)$, the divisor $A = (a,b,c)$ is ample on $S$. Let $C$ be a nonsingular curve on $S$ such that $C \sim A$. Again, $C$ exists as shown in \cite{c}. 

Let $R = \CC[x_1, \ldots, x_4]$ be the coordinate ring of $\PP^3$, and let $I$ be the defining ideal of $C$ in $\PP^3$. We will show that there exist $(a,b,c) \in \ZZ^3$ satisfying (\ref{abc}) and a curve $C$ as above, such that
\begin{equation}\label{eq10}
\lim_{n \to \infty} \sfrac{\l(H^0_\mm(R/I^n))}{n^4} \not\in \QQ, 
\end{equation}
where $\mm = (x_1, \ldots, x_4)$ is the maximal homogeneous ideal of $R$. 
Theorem \ref{Theorem*} is thus an immediate consequence.

As in proved in (\ref{eq3}), (\ref{eq2}) and (\ref{eq4}) of Theorem \ref{main-proof}, 
$$\lim_{n \to \infty} \sfrac{\l(H^0_\mm(R/I^n))}{n^4} = \lim_{n \to \infty}\sfrac{\sigma(n)}{n^4} + \lim_{n \to \infty}\sfrac{\tau(n)}{n^4},$$
where $\lim_{n \to \infty}\sfrac{\tau(n)}{n^4} \in \QQ$. It remains to show that there exist $(a,b,c) \in \ZZ^3$ satisfying (\ref{abc}) and a curve $C$ as above, such that
$$ \lim_{n \to \infty} \sfrac{\sigma(n)}{n^4} \not\in \QQ.$$

Let $\I$ be the ideal sheaf of $I$ on $\PP^3$. Let $\pi: X \to \PP^3$ be the blowing up of $\PP^3$ along the ideal sheaf $\I$. There exists a hyperplane $H'$ of $\PP^3$ such that $H'\cdot S = H$. Let $\H$ be the pull-back to $X$ of $H'$, and $E$ the exceptional divisor of the blowing up. Let $\lambda_1 = a - \sqrt{b^2+c^2}$ and $\lambda_2 = a +\sqrt{b^2 + c^2}$. The following facts were proved in \cite{c}.

\begin{lemma} \label{facts}
Suppose $\lambda_2 > 7$. Then,
\begin{enumerate}
\item $h^0(S, \O_S(mH-nC)) = 0$ if $m < \lambda_2n$.
\item $h^0(S, \O_S(mH-nC)) = \sfrac{1}{2}(mH-nC)^2 + 2$ if $m > \lambda_2n$.
\item $h^1(X, \O_X(m\H-nE)) = 0$ if $m > \lambda_2n$.
\end{enumerate}
\end{lemma}

\begin{proof} (1) follows from  \cite[Remark 6]{c}. (2) is a consequence of \cite[Theorem 7]{c}. (3) follows from \cite[Theorem 9]{c}.  
\end{proof}

It was pointed out in \cite[(11)]{c} that
\begin{align}
& \pi_* \O_X(m\H-nE) \cong \I^n(m), \label{direct-1} \\
& R^i \pi_* \O_X(m\H-nE) = 0, \ \text{for} \ i > 0. \label{direct-2}
\end{align}
Thus, we can use the cohomology groups of $\O_X(m\H-nE)$ to calculate $\sigma(n)$. For convenience, we will use $H^i(X, m\H-nE)$ and $H^i(S, mH-nC)$ to denote $H^i(X, \O_X(m\H-nE))$ and $H^i(S, \O_S(mH-nC))$, respectively. When there is no danger of confusion, we shall further omit the space $X$ and $S$ in these cohomology groups. It was also shown in \cite[(12)]{c} that there exists the following exact sequence:
\begin{align}
0 \to \O_X((m-4)\H-(n-1)E) \to \O_X(m\H-nE) \to \O_S(mH-nC) \to 0. \label{exact-seq}
\end{align}

The existence of the desired example follows from the following theorem.

\begin{theorem} \label{irrationality}
There exist $(a,b,c) \in \ZZ^3$ satisfying (\ref{abc}) and a corresponding nonsingular curve $C$ such that, if $I \subset R$ is the defining ideal of $C$, then
$$ \lim_{n \to \infty} \sfrac{\sigma(n)}{n^4} \not\in \QQ,$$
and
$$
\lim_{n \to \infty} \sfrac{\l(H^0_\mm(R/I^n))}{n^4} \not\in \QQ. 
$$
\end{theorem}

\begin{proof} Taking the long exact sequence of cohomology groups from the exact sequence (\ref{exact-seq}), we get
\begin{align*}
0 \to H^0((m-4)\H-(n-1)E) \to & H^0(m\H-nE) \to \\
& H^0(mH-nC) \to H^1((m-4)\H-(n-1)E).
\end{align*}
It follows from Lemma \ref{facts} that 
\begin{align}
& H^0(m\H-nE) \cong H^0((m-4)\H-(n-1)E) \ \text{if} \ m < \lambda_2n, \label{induction-1} \\
& H^0(m\H-nE) \cong H^0((m-4)\H-(n-1)E) \oplus H^0(mH-nC) \ \text{if} \ m > \lambda_2n. \label{induction-2}
\end{align}

Write $m = 4n+r$. Consider the following cases. \par
\noindent {\it Case 1:} $r < 0$. Since $\lambda_2 > 4$, we have $m < \lambda_2n$. Thus, using (\ref{induction-1}) and successive induction, we get
\begin{align} 
H^0(m\H-nE) = H^0(r\H) = H^0(\O_{\PP^3}(r)) = 0. \label{case1}
\end{align} \par
\noindent {\it Case 2:} $r \ge 0$. If $r > (\lambda_2-4)n$, i.e. $m > \lambda_2n$, then using (\ref{induction-2}) and successive induction, we get
$h^0(m\H-nE) = \sum_{k=1}^n h^0((r+4k)H-kC) + 2n+h^0(\O_{\PP^3}(r)).$
Lemma \ref{facts} now gives
\begin{align} 
h^0(m\H-nE) = \sfrac{1}{2}\sum_{k=1}^n ((r+4k)H-kC)^2 +2n+ h^0(\O_{\PP^3}(r)). \label{case2-1}
\end{align}
On the other hand, if $r < (\lambda_2-4)n$, then put $t = \Big[\sfrac{r}{\lambda_2-4}\Big]$. By successive induction using both (\ref{induction-1}) and (\ref{induction-2}), we get
$$
\begin{array}{ll} 
h^0(m\H-nE) &= \sum_{k=1}^t h^0((r+4k)H-kC) + h^0(\O_{\PP^3}(r))\\
&= \sfrac{1}{2}\sum_{k=1}^t ((r+4k)H-kC)^2 + 2t+h^0(\O_{\PP^3}(r)). \label{case2-2}
\end{array}
$$
By (\ref{direct-1}), we have 
\begin{align*}
\sigma(n) & = \sum_{m=0}^{en} h^0(m\H-nE) \\ 
& = \sum_{r=-4n}^{(e-4)n}h^0((r+4n)\H-nE) \\ 
& = \sum_{r=-4n}^{-1}h^0((r+4n)\H-nE) + \sum_{r=0}^{[(\lambda_2-4)n]}h^0((r+4n)\H-nE) \\
& \quad + \sum_{r=[(\lambda_2-4)n]+1}^{(e-4)n}h^0((r+4n)\H-nE). 
\end{align*}
This together with (\ref{case1}), (\ref{case2-1}) and (\ref{case2-2}) gives us
\begin{eqnarray*}
\sigma(n) & = & \sum_{r=0}^{[(\lambda_2-4)n]} \left(\sfrac{1}{2}\sum_{k=1}^{\big[\frac{r}{\lambda_2-4}\big]}((r+4k)H-kC)^2 
+2\Big[\sfrac{r}{\lambda_2-4}\Big]
+ h^0(\O_{\PP^3}(r))\right) \\
& \quad & + \sum_{r=[(\lambda_2-4)n]+1}^{(e-4)n} \left(\sfrac{1}{2}\sum_{k=1}^n ((r+4k)H-kC)^2 + 2n+h^0(\O_{\PP^3}(r))\right) \\
& = & \sum_{r=0}^{[(\lambda_2-4)n]} \left(\sfrac{1}{2}\sum_{k=1}^{\big[\frac{r}{\lambda_2-4}\big]}((r+4k)H-kC)^2\right) \\
& \quad & + \sum_{r=[(\lambda_2-4)n]+1}^{(e-4)n} \left(\sfrac{1}{2}\sum_{k=1}^n ((r+4k)H-kC)^2\right) + \sum_{r=0}^{(e-4)n} h^0(\O_{\PP^3}(r))\\
& \quad &+2\Big(\sum_{r=0}^{[\lambda_2-4)n]}\Big[\sfrac{r}{\lambda_2-4}\Big]+n((e-4)n-[(\lambda_2-4)n]-1)\Big).
\end{eqnarray*}

 Let 
$$
Q(s,r) = \sfrac{1}{2}\sum_{k=1}^s ((r+4k)H-kC)^2,
$$
$$
V(n)=2\Big(\sum_{r=0}^{[\lambda_2-4)n]}\Big[\sfrac{r}{\lambda_2-4}\Big]+n((e-4)n-[(\lambda_2-4)n]-1)\Big)
$$
 and 
$$
U(n) = \sum_{r=0}^{(e-4)n}h^0(\O_{\PP^3}(r)).
$$
 Let 
$$
P(s,r) = Q(s,r) - Q(s-1,r)=\sfrac{1}{2}((r+4s)H-sC)^2
$$ 
with the convention that $Q(s,r) = 0$ for $s < 1$. For simplicity, let us also denote $\lambda=\lambda_2-4$. Then, we can rewrite $\sigma(n)$ as follows:
\begin{align}
 \sigma(n)=\sum_{r=[\lambda]+1}^{[\lambda n]} P(1,r) + \cdots + \sum_{r=[\lambda (n-1)]+1}^{[\lambda n]} P(n-1,r) + \sum_{r=[\lambda n]+1}^{(e-4)n} Q(n,r) + U(n)+V(n). \label{bigsum}
\end{align}
Let us consider one term $\sum_{r=[\lambda l]+1}^{[\lambda n]} P(l,r)$ of the sum (\ref{bigsum}) for some $1 \le l \le n-1$. From the intersection form $q$ of (\ref{eq4}) we have
\begin{align*} 
\sum_{r=[\lambda l]+1}^{\lambda n]} P(l,r) & = \sum_{r=[\lambda l]+1}^{[\lambda n]} 2\left((r+(4-a)l)^2 - (b^2+c^2)l^2\right) \\
& = \sum_{r=1}^{[\lambda n]} 2\left((r+(4-a)l)^2 - (b^2+c^2)l^2\right) - \sum_{r=1}^{[\lambda l]} 2\left((r+(4-a)l)^2 - (b^2+c^2)l^2\right) \\
& = \left( \sfrac{2}{3}([\lambda n])^3 + 2(4-a)l([\lambda n])^2 + 2((4-a)^2-b^2-c^2)l^2([\lambda n]) \right) \\
& \quad - \left( \sfrac{2}{3}([\lambda l])^3 + 2(4-a)l([\lambda l])^2 + 2((4-a)^2-b^2-c^2)l^2([\lambda l]) \right) \\
& \quad + [\lambda n]^2+(\sfrac{1}{3}+2(4-a)l)[\lambda n]+[\lambda l]^2
+(\frac{1}{3}+2(4-a)l)[\lambda l].
\end{align*}
Note that $\lambda l -1 < [\lambda l] < \lambda l$ for any $l$. Thus, we have
\begin{align*} 
\sum_{r=[\lambda l]+1}^{[\lambda n]} P(l,r) & = \left( \sfrac{2}{3}(\lambda n)^3 + 2(4-a)l(\lambda n)^2 + 2((4-a)^2-b^2-c^2)l^2(\lambda n) \right) \\
& \quad - \left( \sfrac{2}{3}(\lambda l)^3 + 2(4-a)l(\lambda l)^2 + 2((4-a)^2-b^2-c^2)l^2(\lambda l) \right) \\
& \quad + F(n,l) 
\end{align*}
where $F(n,l)$ is a function such that there exists a polynomial $G(n,l)$ of degree 2
with positive real coefficients satisfying $\mid F(n,l)\mid< G(n,l)$ for all $n,l\in{\bf N}$.

Taking the sum as $l$ goes from $1$ to $(n-1)$, we get
\begin{align*}
\sum_{l=1}^{n-1} \sum_{r=[\lambda l]+1}^{[\lambda n]} P(l,r) & = \left( \sfrac{2}{3}\lambda^3 n^4 + (4-a)\lambda^2 n^4 + \sfrac{2}{3}((4-a)^2-b^2-c^2)\lambda n^4\right) \\
& \quad - \left( \sfrac{1}{6}\lambda^3 n^4 + \sfrac{1}{2}(4-a)\lambda^2 n^4 + \sfrac{1}{2}((4-a)^2-b^2-c^2)\lambda n^4\right) + O(n^3).
\end{align*}
We also have
\begin{align*}
Q(n,r) & = \sum_{k=1}^n 2\left((r+(4-a)k)^2 - (b^2+c^2)k^2\right) \\
& = 2r^2n + 2(4-a)rn^2 + \sfrac{2}{3}((4-a)^2-b^2-c^2)n^3 + H(n,r)
\end{align*}
where $H(n,r)$ is a real polynomial of degree $\le 2$ in $n$ and $r$.
Thus,
\begin{align*}
\sum_{r=[\lambda n]+1}^{(e-4)n} Q(n,r) & = \sum_{r=1}^{(e-4)n} Q(n,r) - \sum_{r=1}^{[\lambda n]} Q(n,r) \\
& = \left( \sfrac{2}{3}(e-4)^3n^4 + (4-a)(e-4)^2n^4+\sfrac{2}{3}((4-a)^2-b^2-c^2)(e-4)n^4 \right) \\
& \quad - \left( \sfrac{2}{3}n([\lambda n])^3 + (4-a)n^2([\lambda n])^2+\sfrac{2}{3}((4-a)^2-b^2-c^2)n^3[\lambda n] \right) + O(n^3) \\
& = An^4 - \left( \sfrac{2}{3}\lambda^3 n^4 + (4-a)\lambda^2 n^4+\sfrac{2}{3}((4-a)^2-b^2-c^2)\lambda n^4 \right) + O(n^3),
\end{align*}
where $A \in \QQ$. Moreover, 
$$
U(n) = \sum_{r=0}^{(e-4)n} h^0(\O_{\PP^3}(r)) = \sum_{r=0}^{(e-4)n}
\left(\stackrel{\textstyle r+3}{3}\right)
 = \sfrac{1}{24}(e-4)^4n^4 + O(n^3) = Bn^4 + O(n^3)
$$
 with $B \in \QQ$.
Further, we have $V(n)=\,\, O(n^2)$. Hence,
$$ \sigma(n) = (A+B)n^4 - \left(\sfrac{1}{6}\lambda^3 +\sfrac{1}{2}(4-a)^2\lambda^2 +\sfrac{1}{2}((4-a)^2-b^2-c^2)\lambda\right)n^4 + O(n^3).$$

Finally, take $a=4, b=3$ and $c=2$, then clearly $(a,b,c)$ satisfies all the requirements in (\ref{abc}) and $\lambda_2 = a + \sqrt{b^2+c^2} > 7$. We have $\sigma(n) = (A+B)n^4 + \sfrac{13\sqrt{13}}{3}n^4 + O(n^3)$, where $A,B \in \QQ$. Therefore,
$$ \lim_{n \to \infty} \sfrac{\sigma(n)}{n^4} = A+B+\sfrac{13\sqrt{13}}{3} \not\in \QQ.$$
The theorem is proved.
\end{proof}

We now prove Theorem \ref{examp}. Let $I$ and $R$ be the  ideal and ring of this section.
 Let $\mm=(x_1,\ldots,x_r)$, $S=R_{\mm}$, $J=I_{\mm}$. Let $E$ be the injective hull of ${\CC}$.
$$
\text{Hom}_S(H^0_{\mm S}(S/J^n),E)\cong \ext_S^d(S/J^n,S)
$$
by local duality (c.f. Theorem 3.5.8 \cite{BH}). Thus
$$
\l(\ext_S^d(S/J^n,S))=\l(H^0_{\mm S}(S/J^n))=\l(H^0_{\mm }(R/I^n))
$$
for all $n$ (c.f. Proposition 3.2.12 \cite{BH}).
Now the theorem follows from Theorem \ref{Theorem*}.

%%%%%%%%%%%%%%%%%%%%%%%%%%%%%%%%%%%%%%%%%%%%%%%%%%%%%%%%%%%%%%%%%%%%%%%%%

\end{document}